\documentclass[11pt]{article}
\usepackage[utf8]{inputenc}
\usepackage{latexsym}
\usepackage{amsmath,amsxtra,amssymb,latexsym, amscd,amsthm,}
 \usepackage[pdftex]{color,graphicx}
 \usepackage[all]{xy}
 \usepackage{hyperref}
\usepackage{pgf,tikz}
\usepackage[T1]{fontenc}
\usepackage{cleveref}
\usepackage{amsfonts}
\usetikzlibrary{arrows}
\usepackage{graphicx}
\usepackage{caption}
\usepackage{subcaption}
\usepackage{graphpap}
\usepackage{rotating}
\usepackage{bm}
\theoremstyle{plain}

\usepackage{fancyhdr}
\newcommand\shorttitle{Triangles, Fractales and Spaghetti}
\newcommand\authors{E. A. A. Diop, M. Gaye, A. K. Sane}
\theoremstyle{definition}

\title{\sc{Triangles, Fractales and Spaghetti}}
\author{ElHadji Abdou Aziz Diop\footnote{azizdiop55@gmail.com, UCAD, Senegal}, Masseye Gaye\footnote{masseye.gaye@ucad.edu.sn, UCAD, Senegal}, Abdoul Karim Sane\footnote{abdoulkarim3.sane@ucad.edu.sn, UCAD, Senegal}}
\date{ \small{}}

\begin{document}
\renewcommand{\proofname}{Proof}
\renewcommand{\abstractname}{Abstract}
\renewcommand{\refname}{Bibliography}
\maketitle
\begin{abstract} There is well-known problem of geometric probability which can be quote as the Broken Spaghetti Problem. It addresses the following question: \textit{A stick of spaghetti breaks into three parts and all points of the stick have the same probability to be a breaking point. What is the probability that the three sticks, putting together, form a triangle?}

In these notes, we describe a hidden geometric pattern behind the symmetric version of this problem, namely a fractal that parametrizes the sample space of this problem. Using that fractal, we address the question about the probability to obtain a $\delta$-equilateral triangle.   
\end{abstract}
\begin{flushright}
 \footnotesize{\emph{Les derniers \'echos d'une \'etoile qui s'effondre;\\ Hommage \`a un Recteur}}
\end{flushright}
\begin{section}{Introduction}
The Broken Spaghetti Problem, also called the broken stick problem, is an old mathematical problem. It goes back to British mathematicians and it has interested mathematicians like E. Lemoine \cite{Lem}, H. Poincar\'e \cite{Poinc}, de Finetti ~\cite{Fin1, Fin2}. There are several equivalent ways to formulate this problem but we will keep the one mention above. It happens that the answer is $$\widetilde{\mathbb{P}}=\frac{1}{4};$$ and it was already known to British mathematicians. This problem has been generalized to higher dimension object in \cite{Andr, Eis}, and in \textit{"Things to do with broken sticks"} \cite{Ion} the authors addressed a variant of the Broken Spaghetti Problem by changing triangles with other geometric objects. But, let us stay on the original problem. As mentioned in \cite{Good}, E. Lemoine was the first to publish an article on this problem, and he used an exhaustion process to outcome the answer. There is also a geometric approach given by H. Poincar\'e in \textit{Calcul des Probabilit\'es} \cite{Poinc}, which amounts to describe the space of all possibilities: an equilateral triangle, and then compute the probability as the relative area. The geometrical approach was known before Poincar\'e but what he did better is the proof of why the fact that the broken points are equally likely translate to a uniform distribution on the sample space. This  allows to compute the probability as relative area and it was assume to be evident by other authors. In \cite{Good}, G. S. Goodman gave another proof of why the samples are uniformly distributed using a beautiful argument from elementary geometry. Since Lemoine method uses uniform distribution on a discretized version of the problem and a limit process, his result combine with Poincar\'e's one could be interpreted like a convergence of a discrete uniform distribution to a continuous one.        

Goodman also raised the question on the choice on how to sample this problem and how it affects the outcome of the probability. The problem is mostly sample as followed (see Figure \ref{samp}):
\begin{itemize}
\item $l_1$, $l_2$ and $l_3$ are the length of the three sticks
\item $l_1$ $l_2$ and $l_3$ are ordered in such a way that $l_1$ is the length of stick on the left, $l_2$ for the middle sitck and $l_3$ for the one on the right. 
\end{itemize}
\begin{figure}[htbp]
\begin{center}
\includegraphics[scale=0.4]{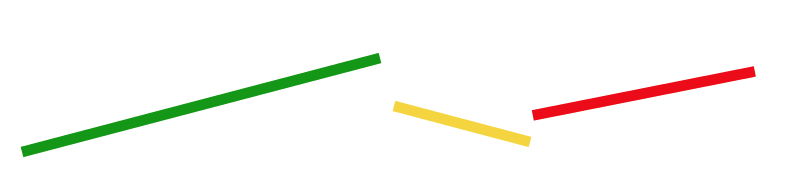}
\put(-40,3){$l_1$} \put(-25,2){$l_2$} \put(-10,4){$l_3$}
\caption{Sampling of the (non symmetric) spaghetti problem}
\label{samp}
\end{center}
\end{figure}

Using this way, the sample $l_1=3/4, l_2=1/8, l_3=1/8$ is considered to be different to $l_1=1/8, l_2=3/4, l_3=1/8$. Although, one can consider a symmetric way to sample this problem and in this case, the order on which the three sticks appears is worthless. Let us call this problem the symmetric Broken Spaghetti Problem.

\textit{In these notes, we describe the geometrical shape of the sample space for the symmetric Broken Spaghetti Problem. It happens that it is a fractal and the probability for the symmetric version is:
$$\mathbb{P}=\frac{1}{8}.$$ The interpretation of the fractal enable to give the probability to obtain a $\delta$-equilateral triangle.}  

Nowadays, the Broken Spaghetti Problem is also used to introduce the notion of probability to high school students. The authors have presented this problem during the BRIS-NLAGA 2022 held in Ziguinchor/Senegal. The goal was to introduce the notion of probability to students and the notion of fractal as well.  

As one can see, the symmetric version of the Broken Spaghetti Problem does note have the same probability like the original problem. This brings to mind the Bertrand paradox in which a problem have different probability depending on the way the sampling is made. Nonetheless, our goal was not to insist on the paradoxical behavior of the Broken Spaghetti Problem but rather to attach a fractal to this problem and thereby to show to students how rich this simple problem is.  By the way, we would like to end up this introduction with a beautiful  sentence of G. S. Goodman \cite{Good} on this problem:
\begin{center}
\textit{\textbf{...Because of this, the problem of the Broken Stick, often snubbed as a mere mathematical diversion by those who forget that probability theory had its origins in mathematical diversions, deserves to occupy a more dignified place in the hierarchy of mathematical though.}}
\end{center}
\end{section}
\begin{section}{Geometric approach reviewed}
In this section, we set up some materials and we give the geometric proof of the Broken Spaghetti Problem. Assume that the stick has length $1$ and it breaks randomly at two points. Let $l_1$, $l_2$ and $l_3$ be the (ordered) list of length of the three sticks. Therefore, the following equation holds: 
\[l_1+l_2+l_3=1\quad (1).\]
So, a sample here is a vector $v:=(l_1, l_2, l_3)\in\mathbb{R}^3$ such that $(1)$ is satisfied. It follows that the set of all possibilities
$$\widetilde{\mathcal{E}}=\{(l_1,l_2,l_3)\in \mathbb{R}^*_+\times\mathbb{R}^*_+\times\mathbb{R}^*_+, l_1+l_2+l_3=1\},$$
is an equilateral triangle: the two dimensional simplex in $\mathbb{R}^3$.  There is an another way to represent $\widetilde{\mathcal{T}}$ just by drawing it in $\mathbb{R}^2$. In this case, $\widetilde{\mathcal{T}}$ is an equilateral triangle in $\mathbb{R}^2$ with side-length equal to $\frac{2\sqrt{3}}{3}$ and the coordinates of a point $M$ in  $\widetilde{\mathcal{T}}$  is given by the distances between $M$ and each of the three sides (Figure \ref{ssp}).

\begin{figure}[htbp]
\begin{center}
\includegraphics[scale=0.22]{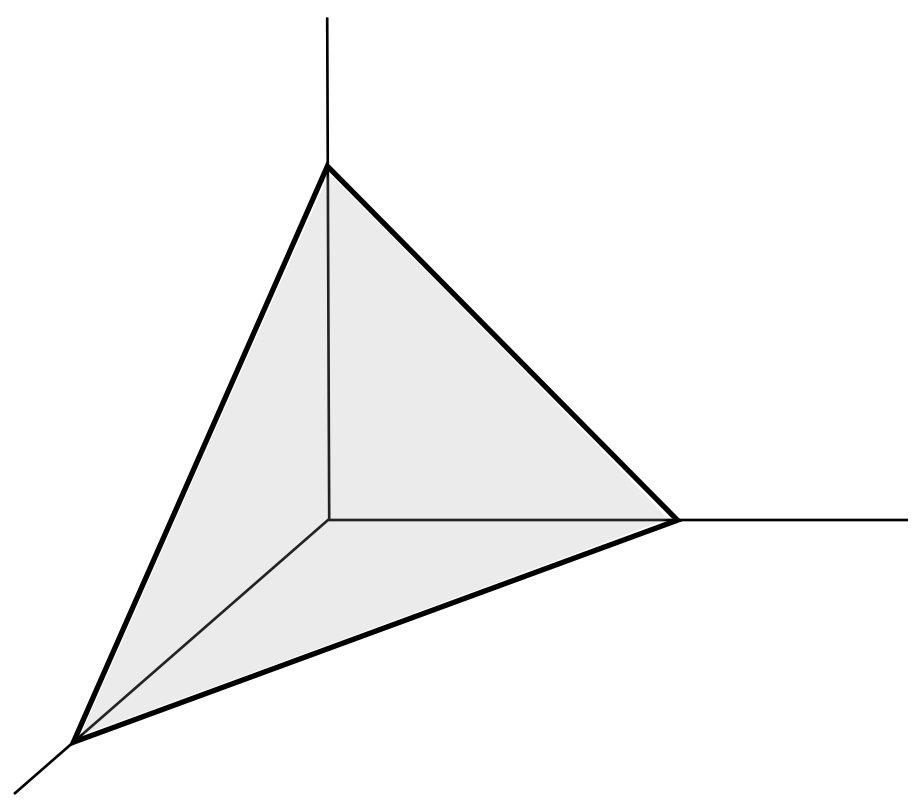}\hspace{1cm}\includegraphics[scale=0.22]{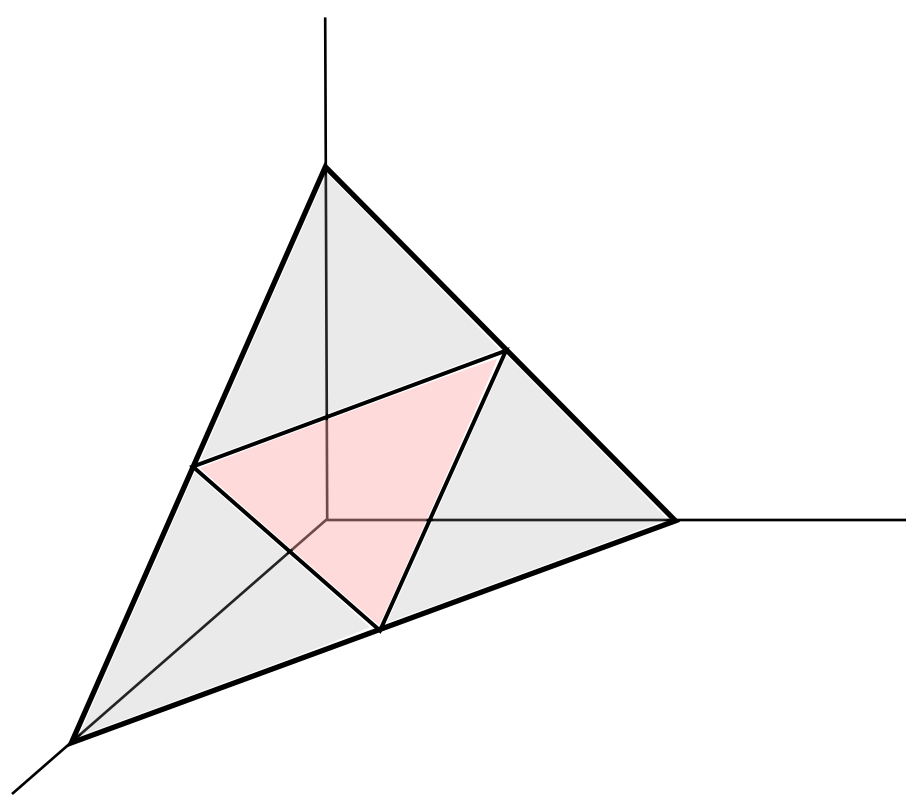}\hspace{1cm}\includegraphics[scale=0.3]{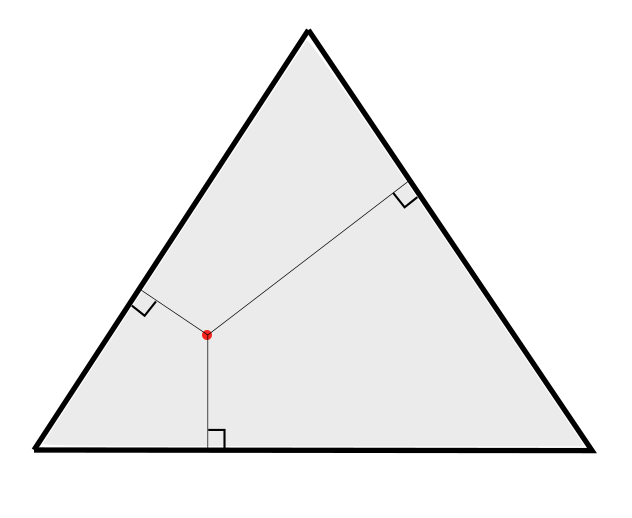}\put(-25, 12){\tiny{$l_1$}} \put(-15, 14){\tiny{$l_3$}} \put(-22, 6){\tiny{$l_2$}}
\caption{The sample space and its representation in $\mathbb{R}^2$.}
\label{ssp}
\end{center}
\end{figure}

Moreover $l_1$, $l_2$ and $l_3$ are length of a triangle providing they satisfy the triangle inequalities: 
\[l_2+l_3\leq l_1,\quad l_1+l_3\leq l_2,\quad l_1+l_2\leq l_3\quad (2).\]
These conditions are equivalent to the following:
\[l_1\leq \frac{1}{2},\quad l_2\leq \frac{1}{2},\quad l_3\leq \frac{1}{2}.\] 

Then, the sample that give a positive answer to the problem are: 
$$\widetilde{\mathcal{T}}=\{(l_1,l_2,l_3)\in \mathbb{R}^*_+\times\mathbb{R}^*_+\times\mathbb{R}^*_+; l_1\leq\frac{1}{2},\hspace{0.2cm} l_2\leq\frac{1}{2},\hspace{0.2cm} l_3\leq\frac{1}{2}\}.$$

The set $\widetilde{\mathcal{T}}$ is also an equilateral triangle: the one that joins the middle of the side of $\widetilde{\mathcal{E}}$. 

By Poincar\'e \cite{Poinc}, the fact that the sticks broke randomly at two points translate to a uniform distribution $\widetilde{\mathcal{E}}$. So, the probability $\widetilde{\mathbb{P}}$ is the amount of $\widetilde{\mathcal{T}}$ we have in $\widetilde{\mathcal{E}}$ (relative area):

\[\widetilde{\mathbb{P}}=\displaystyle{\frac{\mathrm{Area}(\widetilde{\mathcal{T}})}{\mathrm{Area}(\widetilde{\mathcal{E}})}}=\frac{1}{4}.\]   

Now, let us turn to the symmetric spaghetti problem. For that, we consider samples $(l_1, l_2, l_3)$ up to the action of the symmetric group $\mathcal{S}_3$:
\[\sigma.(l_1, l_2, l_3)=(l_{\sigma(1)}, l_{\sigma(2)}, l_{\sigma(3)});\]
where $\sigma\in\mathcal{S}_3$.  

Therefore, the sample space $\mathcal{E}$ of the symmetric spaghetti problem is:
\[\mathcal{E}=\widetilde{\mathcal{E}}/\mathcal{S}_3.\]

\end{section}
\begin{section}{The geometry of the sample space}
In this section, we describe the sample space $\mathcal{E}$. Since $\mathcal{S}_3$ is generated by $(1 2)$, $(1 3)$ and $(2 3)$, the action of $\mathcal{S}_3$ on $\widetilde{\mathcal{E}}$ allow us to consider the points up to symmetries along the three bisectors of $\widetilde{\mathcal{E}}$. Using this, we describe $\mathcal{E}$ inductively.  

\begin{figure}[htbp]
\begin{center}
\includegraphics[scale=0.25]{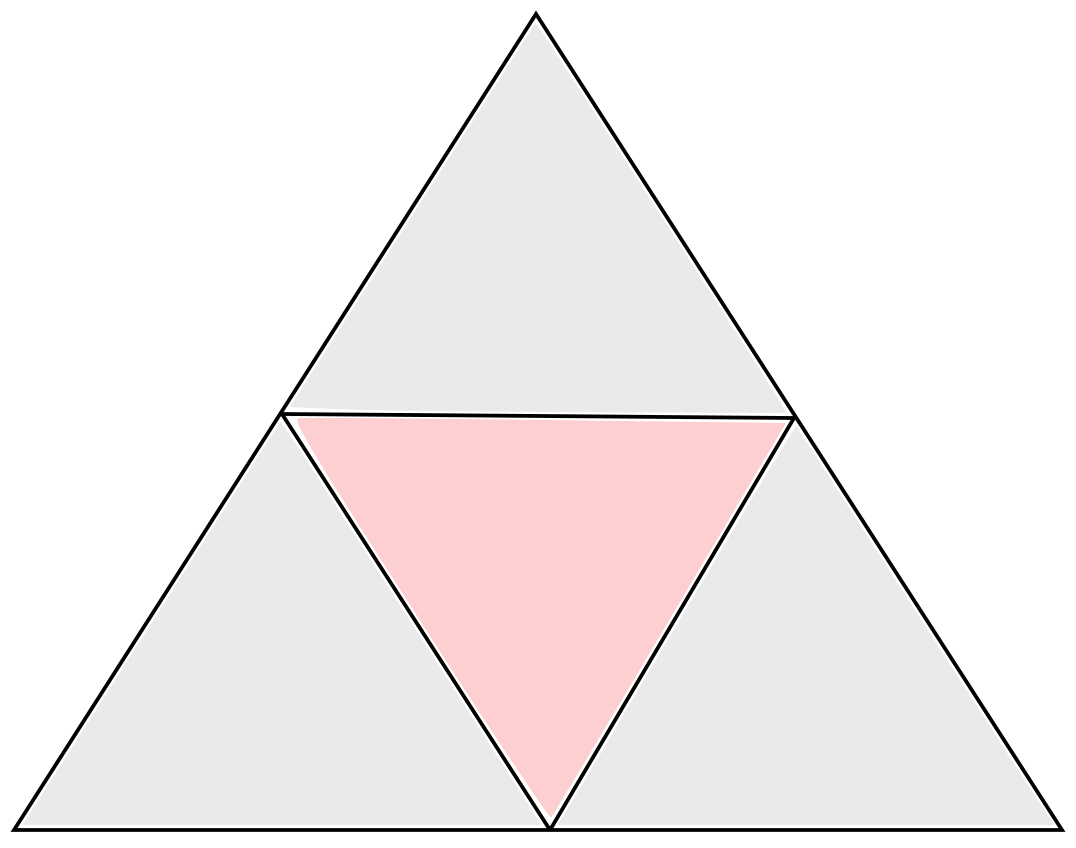}\hspace{1,5cm}\includegraphics[scale=0.25]{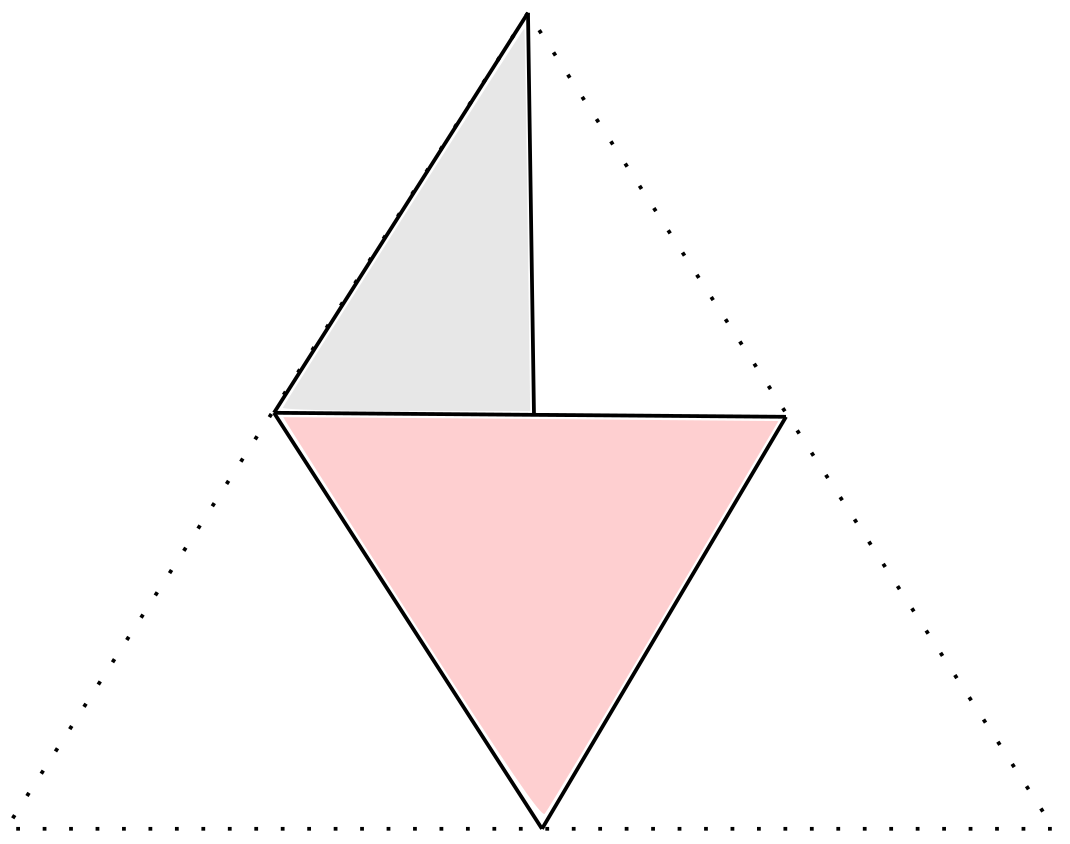}
\put(-60,15){\Large{$\longrightarrow$}}
\caption{Sample space after one step}
\label{1}
\end{center}
\end{figure}
\begin{paragraph}{Step 1:} Let $A_1$ be the triangle joining the middle of the side of  $\widetilde{\mathcal{E}}$. Then,  $\widetilde{\mathcal{E}}$ is separate in three triangles $T_1$, $T_2$ and $T_3$ as depicted in Figure \ref{1}. The triangle $T_i$ is the set of points $(l_1, l_2, l_3)\in\widetilde{\mathcal{E}}$ such that $l_i\geq 1/2$. So a point $A_i$ is equivalent to a point in $A_j$; that is we can delete two of the three triangles (let us say $T_1$ and $T_2$) from $\widetilde{\mathcal{E}}$. In $T_3$, a point $(l_1, l_2, l_3)$ is equivalent to $(l_2, l_1, l_3)$. So, we can delete half the triangle $T_3$. At this step, we obtain a space $\widetilde{\mathcal{E}}_1$ made with two triangles one of which is free from relations (see Figure \ref{1}).  
\end{paragraph}  
\begin{figure}[htbp]
\begin{center}
\includegraphics[scale=0.25]{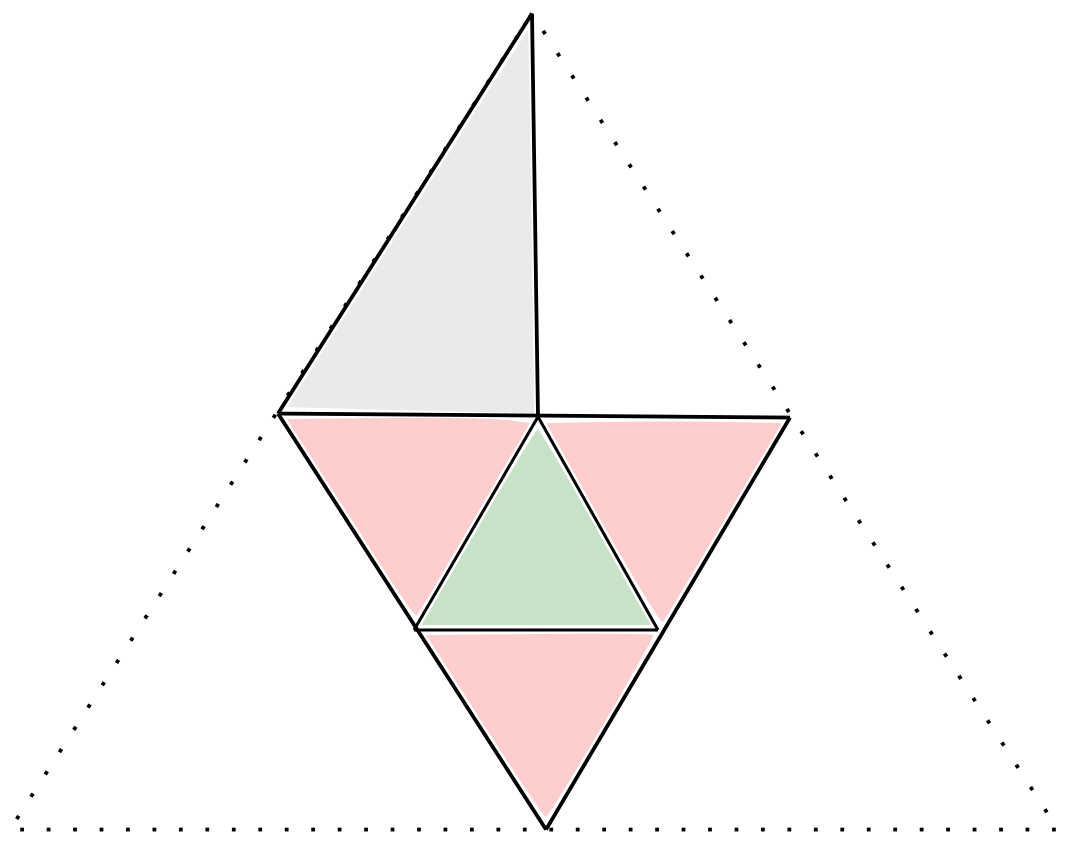}\hspace{1,5cm}\includegraphics[scale=0.25]{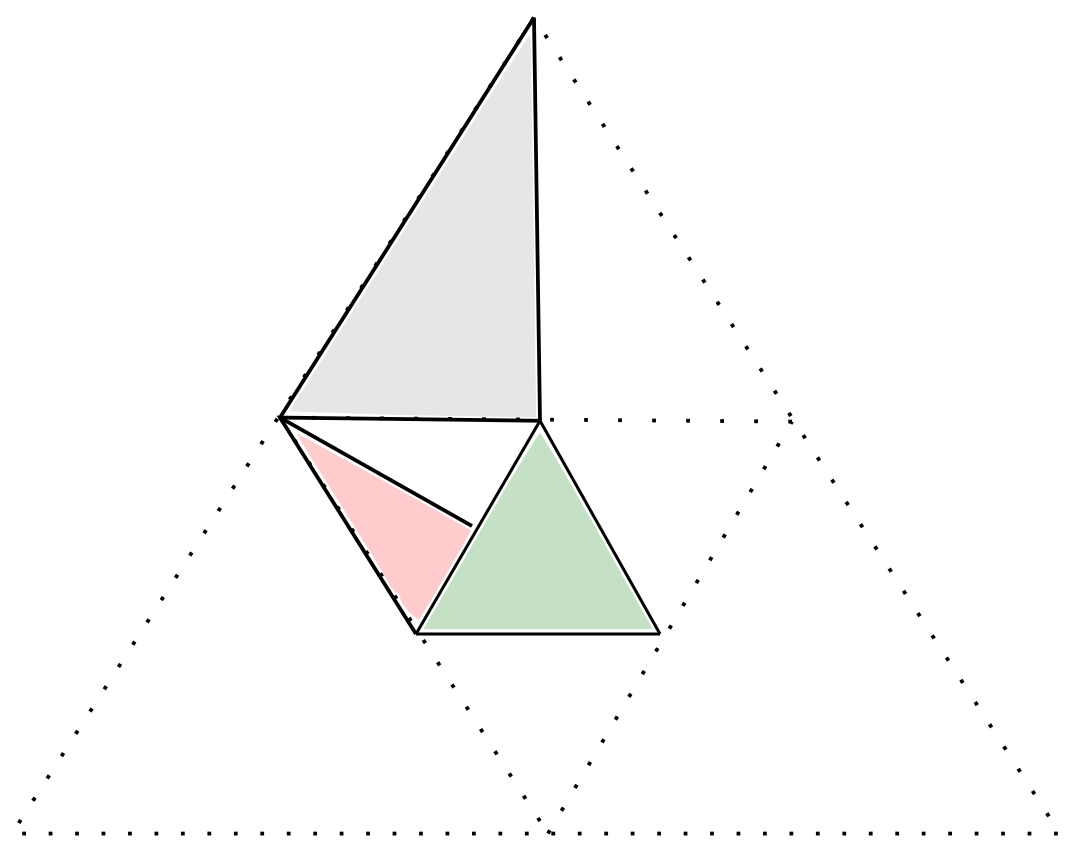}
\put(-60,15){\Large{$\longrightarrow$}}
\caption{Sample space after two steps}
\label{2}
\end{center}
\end{figure} 
\begin{paragraph}{Step 2:}  From $\widetilde{\mathcal{E}}_1$, let $A_2$ be the triangle joining the middle of the sides of $A_1$. Then, $A_2$ split into four triangles $A_2$, $T_1$, $T_2$ and $T_3$ as well. Again, by applying symmetries, we have $T_1\sim T_2\sim T_3$. So, one can remove two of the three triangles and also half of the last one to obtained $\widetilde{\mathcal{E}}_2$ (see Figure \ref{2}).   
\end{paragraph}

\begin{paragraph}{Step n:} At this step, we divide $A_n$ in to four triangles by drawing $A_{n+1}$, the triangle which connects the middle of the sides of $A_n$. Then, $A_n=A_{n+1}\cup T_1\cup T_2\cup T_3$ and $T_1\sim T_2\sim T_3$. We construct $\widetilde{\mathcal{E}}_n$ by deleting $T_2$ and $T_3$ and by removing half of $T_1$.  

\begin{figure}[htbp]
\begin{center}
\includegraphics[scale=0.25]{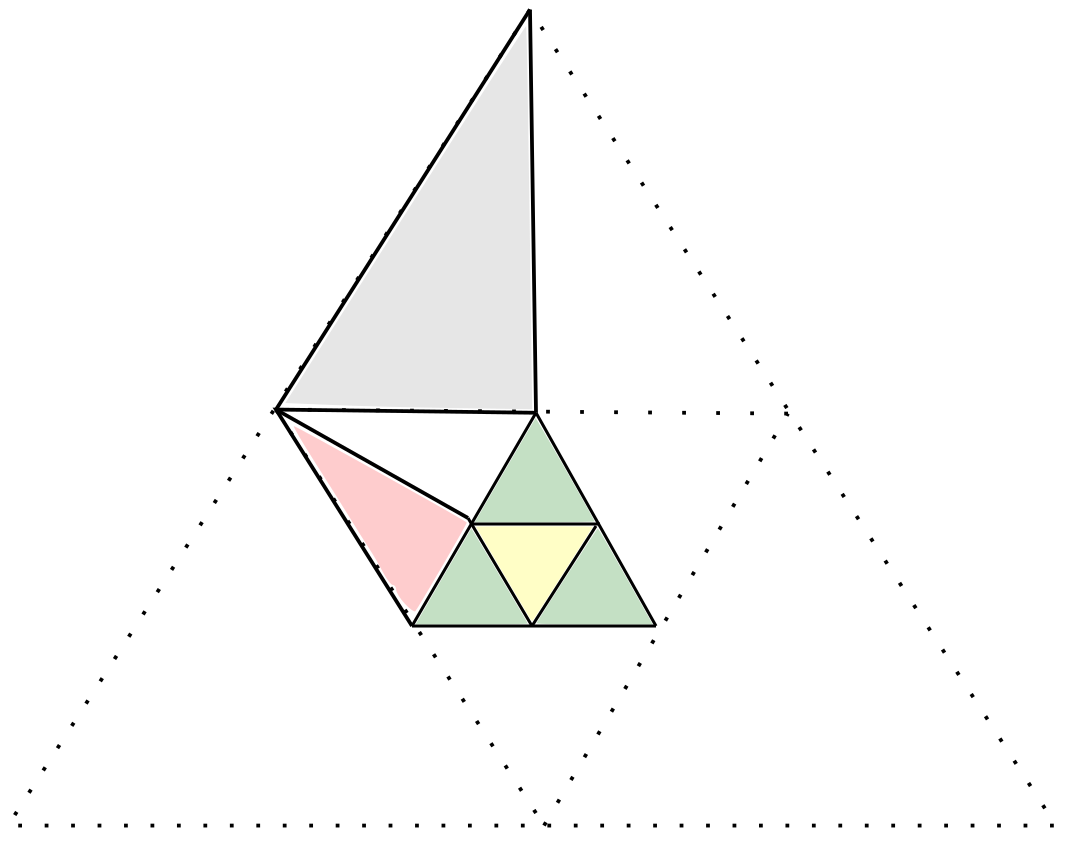}\hspace{1,5cm}\includegraphics[scale=0.25]{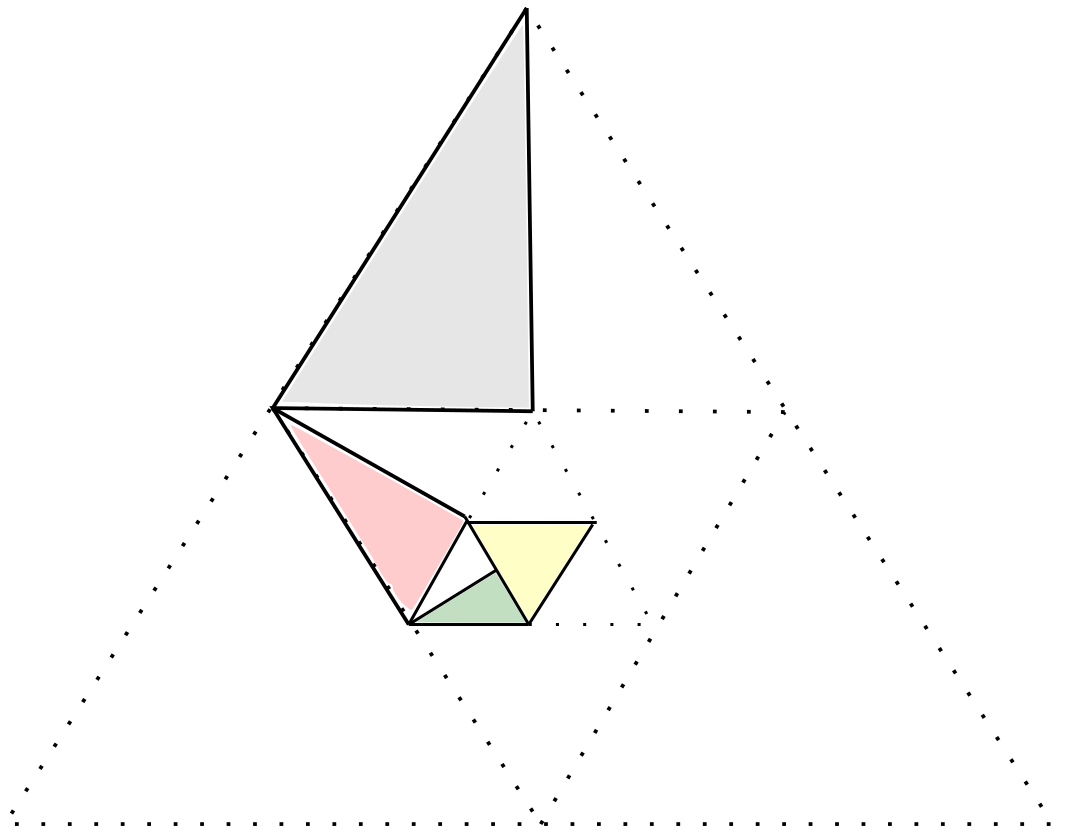}
\put(-60,15){\Large{$\longrightarrow$}}
\caption{Sample space after three steps}
\label{default}
\end{center}
\end{figure}
\end{paragraph}

The Sample space is then given by:
$$\mathcal{E}:=\underset{\underset{n}{\longrightarrow}}{\lim}\hspace{0.1cm}\widetilde{\mathcal{E}}_n.$$

Depending on the pieces one delete during the induction, we obtain different shape of the sample space. In Figure \ref{lim}, we draw two examples of the sample space. 
\begin{figure}[htbp]
\begin{center}
\includegraphics[scale=0.3]{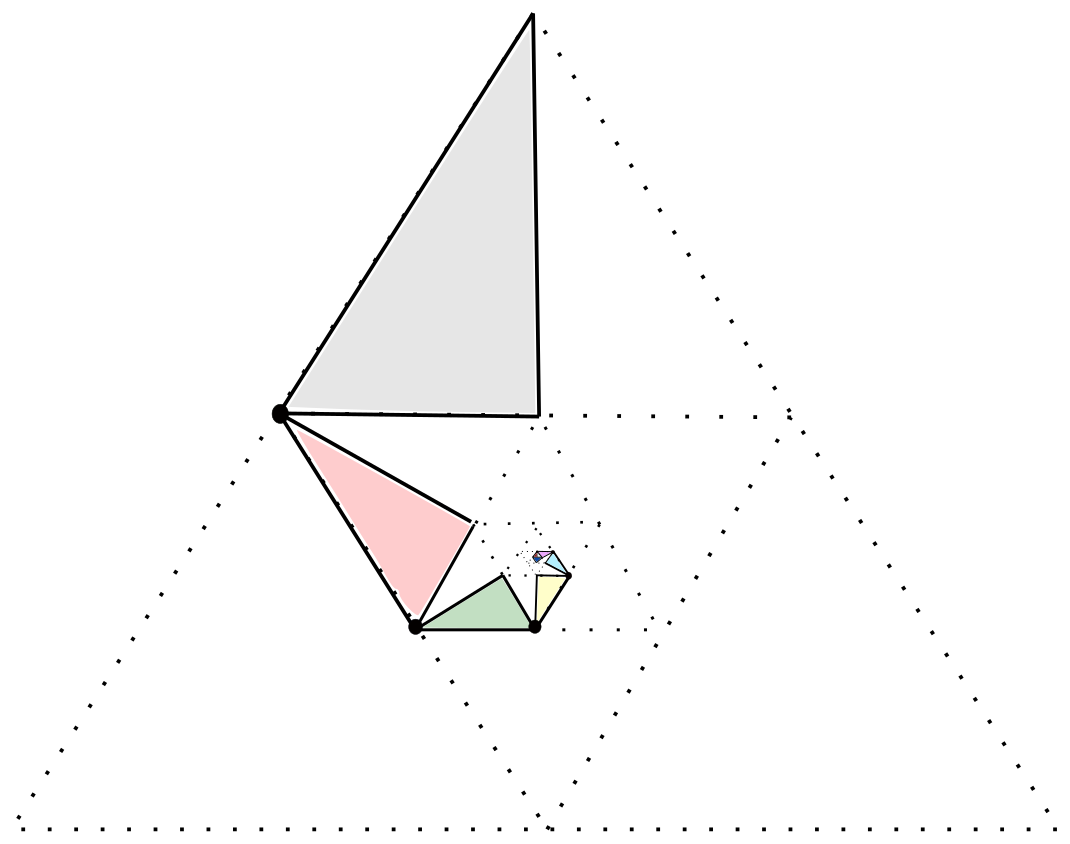} \hspace{1cm} \includegraphics[scale=0.3]{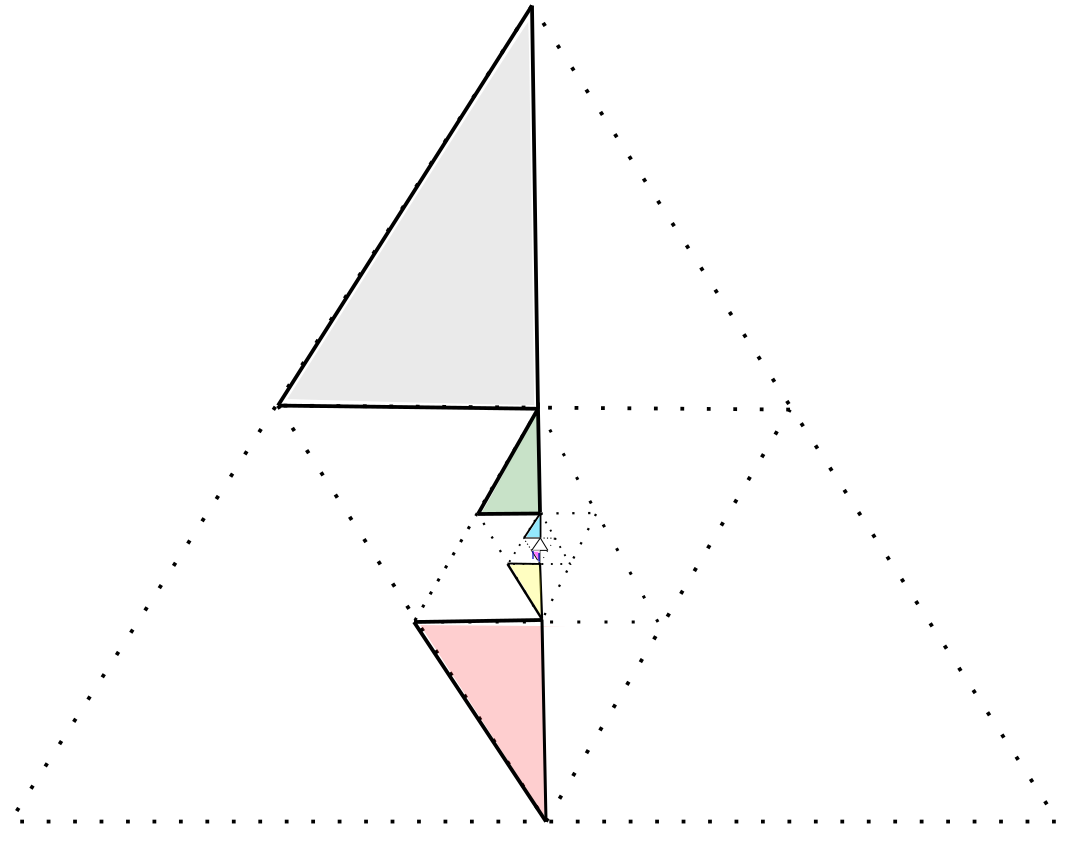}
\caption{Sample space of the symmetric Broken Spaghetti Problem represented in two different ways}
\label{lim}
\end{center}
\end{figure}
\end{section}
\begin{section}{Interpretation of the new sample space}
Since the sample space changed, a natural question is to know if the probability of the non symmetric problem changed as well. Actually, it changes as we will see right away. The sample space $\mathcal{E}$ is the union of infinitely many triangles  $(T_i)_{i=1,\dots,\infty}$. Then: $$\mathrm{Area}(\mathcal{E})=\displaystyle{\sum_{i=1}^\infty \mathrm{Area}(T_i)}.$$ 

The first triangle is one-eight of $\widetilde{\mathcal{E}}$, and each triangle $T_i$ is also one-eight of $T_{i-1}$. Therefore, we have $\mathrm{Area}(T_n)=(\frac{1}{8})^n\mathrm{Area}(\widetilde{\mathcal{E}})$. It follows that:
$$\mathrm{Area}(\mathcal{E})=(\frac{1}{8}+\frac{1}{8^2}+\dots+\frac{1}{8^n}+\dots)\mathrm{Area}(\widetilde{\mathcal{E}})=\frac{1}{7}\mathrm{Area}(\widetilde{\mathcal{E}})$$

The area of the set $\mathcal{T}$ of samples that satisfies the three triangles inequalities is:
$$\mathrm{Area}(\mathcal{T})=(\frac{1}{8^2}+\dots+\frac{1}{8^n}+\dots)\mathrm{Area}(\widetilde{\mathcal{E}})=\frac{1}{7\times8}\mathrm{Area}(\widetilde{\mathcal{E}})$$

The probability to obtain a triangle is now: $$\mathbb{P}=\frac{1}{8}.$$

As one can see, $\mathcal{E}$ is a sequence of triangle pieces that converges to the point $$G:=(\frac{1}{3}, \frac{1}{3}, \frac{1}{3}).$$

This agree with what one expected to have that is the probability to obtain an equilateral triangle is $0$. Let us give the interpretation  of triangles pieces in $\mathcal{E}$. A triangle $\tau$ given by $(l_1, l_2, l_3)$ is said \textit{$\delta$-equilateral} (respectively \textit{$(\delta,\delta')$-equilateral}) if $\max\{|l_1-l_2|, |l_1-l_3|, |l_2-l_3|\}\leq \delta$ (respectively $\delta\leq\max\{|l_1-l_2|, |l_1-l_3|, |l_2-l_3|\}\leq\delta'$).

For each triangle $T_i$, set: $$\delta_i:=\sup\{\max\{|l_1-l_2|, |l_1-l_3|, |l_2-l_3|\}, (l_1, l_2, l_3)\in T_i\}.$$

Then, each triangle piece  $T_i$ ($i\geq2$) in $\mathcal{E}$ represents the sample of points that give a $(\delta_{i+1},\delta_i)$-equilateral triangle and the truncated sequence starting from a piece $T_i$ represents the sample of points that give a $\delta_i$-equilateral triangle.

So, 
\[\mathbb{P}_{i+1,i}=\frac{1}{7\times8^i}, \quad \mathbb{P}_i=\frac{1}{8^{i-1}}; \]    
where $\mathbb{P}_{i+1,i}$ and $\mathbb{P}_i$ are the probability to obtained a $(\delta_{i+1}, \delta_i)$-equilateral triangle and a $\delta_i$-equilateral triangle, respectively. 

\begin{paragraph}{Acknowledgements:} The authors discovered this old problem when one of their student: Ila Wague, asked them this question while he was preparing the annual Fastef-UCAD exam. The authors would like to thank Ila Wague for his question. Also, the authors would like to thank Pr Diaraf Seck for giving them the opportunity to participate to the BRIS-NLAGA 2022 and explain this problem to high school students.  
\end{paragraph}
\end{section}

\end{document}